\documentclass[twoside, 12pt]{article}
\usepackage{amssymb, amsmath, mathrsfs, amsthm}
\usepackage{graphicx}
\usepackage{color}
\usepackage[top=1in, bottom=1in, left=1in, right=1in]{geometry}
\usepackage{float, caption, subcaption}
\usepackage{amsfonts,amsmath,amssymb,amsthm}
\usepackage{graphicx}
\usepackage[linesnumbered,ruled,vlined]{algorithm2e} 
\usepackage{fixmath}
\usepackage{hyperref}
\usepackage{tikz}
\usepackage{makecell} 

\newcommand{\old}[1]{{}}

\def\emptyset{\varnothing}

\usepackage{array}
\newcolumntype{L}[1]{>{\raggedright\let\newline\\\arraybackslash\hspace{0pt}}b{#1}}
\newcolumntype{C}[1]{>{\centering\let\newline\\\arraybackslash\hspace{0pt}}b{#1}}
\newcolumntype{R}[1]{>{\raggedleft\let\newline\\\arraybackslash\hspace{0pt}}m{#1}}

\providecommand{\keywords}[1]{\textbf{\textit{Keywords---}} #1} 
\newcolumntype{M}[1]{>{\centering\arraybackslash}m{#1}}
\begin{document}

\pagestyle{plain}

\title{An efficient algorithm to test forcibly-biconnectedness of graphical degree sequences}

\author{
Kai Wang\footnote{Department of Computer Sciences,
Georgia Southern University,
Statesboro, GA 30460, USA 
\tt{kwang@georgiasouthern.edu}}
}
\maketitle

\begin{abstract}
We present an algorithm to test whether a given graphical degree sequence is forcibly biconnected or not
and prove its correctness. The worst case run time complexity of the algorithm is shown to be exponential but still
much better than the previous basic algorithm presented in \cite{Wang2018}. We show through experimental evaluations
that the algorithm is efficient on average. We also adapt Ruskey et al's classic algorithm to enumerate zero-free
graphical degree sequences of length $n$ and Barnes and Savage's classic algorithm
to enumerate graphical partitions of an even integer $n$ by incorporating our testing algorithm into theirs and then obtain some
enumerative results about forcibly biconnected graphical degree sequences of given length $n$ and forcibly biconnected graphical
partitions of given even integer $n$. Based on these enumerative results we make some conjectures such as:
when $n$ is large, (1) the proportion of forcibly biconnected graphical degree sequences of length $n$ among all zero-free graphical
degree sequences of length $n$ is asymptotically a constant between 0 and 1; (2) the proportion of forcibly biconnected graphical
partitions of even $n$ among all forcibly connected graphical partitions of $n$ is asymptotically 0.

\end{abstract}

\section{Introduction}
We consider graphical degree sequences of finite simple graphs (i.e. finite undirected graphs without loops or multiple edges)
where the order of the terms in the sequence does not matter. As such, the terms in a graphical degree sequence is often written
in non-increasing order for convenience.
An arbitrary non-increasing sequence of non-negative integers $a_1\ge a_2 \ge \cdots \ge a_n$ can be easily tested
whether it is a graphical degree sequence by using the Erd{\H{o}}s-Gallai criterion \cite{ErdosCallai1960} or
the Havel-Hakimi algorithm \cite{Havel1955,Hakimi1962}. Sierksma and Hoogeveen \cite{Sierksma1991} later summarized
seven equivalent criteria to characterize graphical degree sequences.
A zero-free graphical degree sequence is also called a graphical partition. The former terminology is often used when
the length of the sequences under consideration is fixed while the later is often used when the sum of the terms in the partitions
under consideration is fixed.

Enumerating all the graphs having the same vertex degree sequence and exploring their properties have been of interest.
To our knowledge, an efficient algorithm of Meringer \cite{Meringer1999} is available for the case of regular graphs,
but no efficient algorithm is known to solve the problem for general graphical degree sequences. When considering all realizations of
the same vertex degree sequence, two notions are very useful. Let P be any property of graphs (e.g. biconnected, critical,
Hamiltonian, etc). A graphical degree sequence $\mathbf{d}=(d_1\ge d_2 \ge \cdots \ge d_n)$ is called \textit{potentially}
P-graphic if it has at least one realization having the property P and \textit{forcibly} P-graphic if all its realizations have the property P
\cite{Rao1981}. In a previous paper \cite{Wang2018} we have presented an efficient algorithm to test whether a graphical degree sequence is
forcibly connected or not and also outline an algorithmic framework to test whether a graphical degree sequence is forcibly $k$-connected or not for every fixed $k\ge 2$. In this paper we will present a more sophisticated and efficient
algorithm to test forcibly biconnectedness of graphical
degree sequences. Recall that Wang and Cleitman \cite{WangKleitman1973} have given a simple characterization of potentially
$k$-connected graphical degree sequences of length $n$, which can be performed in $O(n)$ time given any input. However,
testing forcibly $k$-connectedness appears to be much harder.
Some sufficient (but unnecessary) conditions are known for a graphical degree sequence
to be forcibly connected or forcibly $k$-connected \cite{chartrand_kapoor_kronk_1968,BOESCH1974,CHOUDUM1991}.

In the rest of this paper we review the basic algorithm and present the improved algorithm to characterize forcibly biconnected
graphical degree sequences and give a proof why it works in Section \ref{sec:alg}. We analyze the complexity of the algorithm
in Section \ref{sec:analysis}. In Section \ref{sec:experiments} we demonstrate the efficiency of the algorithm through
some computational experiments and then present some enumerative results regarding forcibly biconnected graphical degree
sequences of given length $n$ and forcibly biconnected graphical partitions of given even integer $n$.
Based on the observations on these enumerative results we make some conjectures about the relative asymptotic behavior
of considered functions and the unimodality of certain associated integer sequences in Section \ref{sec:conjectures}.
Finally we conclude in Section \ref{sec:conclusion}.

\section{The testing algorithm}
\label{sec:alg}
In this section we first review the basic algorithm to test forcibly biconnectedness proposed in \cite{Wang2018} and then
present an improved version of the algorithm. We will give a detailed proof why it correctly identifies forcibly biconnected graphical
degree sequences and also comment on certain implementation issues.
\subsection{Review of the basic algorithm}
The basic algorithm from \cite{Wang2018} to test forcibly biconnectedness of graphical degree sequences is shown in Algorithm
\ref{algf2_basic}. The
idea is simple as follows. First we need to make sure that the input graphical degree sequence $\mathbf{d}$ is potentially biconnected
and forcibly connected before we continue the test of forcibly biconnectedness. This is why we have lines 1 and 2.

After we confirmed this is the case, we can conclude that the input $\mathbf{d}$ is forcibly biconnected as long as we
become sure that in every realization of $\mathbf{d}$
there is no cut vertex based on the definition of biconnected (non-separable) graphs. In order to simplify our discussion, we call the
degree of a cut vertex in a connected graph a \textit{cut degree}. Recall that we have defined in \cite{Wang2018} that a generalized
Havel-Hakimi (GHH) operation on a graphical degree sequence $\mathbf{d}$ means selecting an arbitrary term $d_i$ to remove
from $\mathbf{d}$ and then selecting an arbitrary collection $\mathbf{d_S}$ of size $d_i$ from the remaining sequence
$\mathbf{d}-\{d_i\}$ and decrementing each of them by 1. The resulting sequence is notationally written as $\mathbf{d'}=
GHH(\mathbf{d},d_i,\mathbf{d_S})$. If $\mathbf{d'}$ is a non forcibly connected graphical degree sequence under some choice of
$d_i$ and $\mathbf{d_S}$, then the $d_i$ is a cut degree and $\mathbf{d}$ is not forcibly biconnected. This is because if $\mathbf{d'}$
is non forcibly connected, i.e. it has a disconnected realization, then $\mathbf{d}$ has a connected realization with a cut vertex of degree
$d_i$, which shows $\mathbf{d}$ is not forcibly biconnected. If no such choice of $d_i$ and $\mathbf{d_S}$ exists, then no cut degree
exists and $\mathbf{d}$ is forcibly biconnected. As indicated in \cite{Wang2018} if $\mathbf{d'}$ obtained on line 4 is not a graphical
degree sequence, then the \textbf{for} loop continues to iterate without returning \textit{False} on line 6.

\begin{algorithm}[h]
	\KwIn{A zero-free graphical degree sequence $\mathbf{d}=(d_1\ge d_2 \ge \cdots \ge d_n)$}
	\KwOut{\textit{True} or \textit{False}, indicating whether $\mathbf{d}$ is forcibly biconnected or not}
	\uIf{$\mathbf{d}$ is not potentially biconnected or forcibly connected}{
		\Return{\textit{False}}
	}
	\For{each $d_i$ and each collection $\mathbf{d_S}$ of size $d_i$ from $\mathbf{d}-\{d_i\}$} {
		$\mathbf{d'} \gets GHH(\mathbf{d},d_i,\mathbf{d_S})$\;
		\uIf{$\mathbf{d'}$ is a non forcibly connected graphical degree sequence}{
			\Return{\textit{False}}
		}
	}
	\Return{\textit{True}}
	\caption{basic version pseudo-code to test whether a graphical degree sequence is forcibly biconnected (See text for the description of GHH operation)}
	\label{algf2_basic}
\end{algorithm}

We remark that although Algorithm \ref{algf2_basic} works, its performance is probably poor due to the large number of possible
choices of $d_i$ and $\mathbf{d_S}$ in the \textbf{for} loop from lines 3 to 6. See Section \ref{sec:analysis} for more detailed worst
case run time complexity analysis. Our rudimentary implementation of Algorithm \ref{algf2_basic}
can start to encounter bottlenecks for input sequences $\mathbf{d}$ of length 30 to 40.

\subsection{Pseudo-code of the improved version and the proof of its correctness}
\label{subsec:pseudo-code}
Our strategy to improve Algorithm \ref{algf2_basic} is to replace the \textbf{for} loop from lines 3 to 6 with more sophisticated
methods to find a potential cut degree. The idea is still simple. Assume there is a cut degree in the input sequence $\mathbf{d}$. Then
we can construct two graphical degree sequences from $\mathbf{d}$ with total length $n+1$ of which $n-1$ come directly from $\mathbf{d}$
and the remaining two constitute a partition of the cut degree. We actually first test some necessary conditions for the existence of
a cut degree, during which
process an auxiliary set $\mathbf{S}$ of the orders of potential smaller induced subgraphs (explained below) is constructed should
a cut degree exist. The pseudo code of the improved version is shown in Algorithm \ref{algf2_improved}.

\begin{algorithm}[!h]
	\KwIn{A zero-free graphical degree sequence $\mathbf{d}=(d_1\ge d_2 \ge \cdots \ge d_n)$}
	\KwOut{\textit{True} or \textit{False}, indicating whether $\mathbf{d}$ is forcibly biconnected or not}
	\uIf{$\mathbf{d}$ is not potentially biconnected or forcibly connected}{
		\Return{\textit{False}}
	}
	\uIf{$d_2+d_n\ge n$}{
		\Return{\textit{True}}
	}
	$\mathbf{S} \gets \emptyset$\; 
	\For{$s \gets d_n$ \textbf{to} $\lfloor (n-1)/2 \rfloor$} {
		\uIf{$d_{n-s+1}\le s$ \textbf{and} $d_2\le n-s-1$}{
			$\mathbf{S} \gets \mathbf{S}\cup \{s\}$
		}
	}
	\uIf{$\mathbf{S}=\emptyset$}{
		\Return{\textit{True}}
	}
	\ForEach{distinct $d$ in $\mathbf{d}$}{
		\ForEach{$s$ in $\mathbf{S}$}{
			\uIf{$d=d_1$ \textbf{or} $n-s-1\ge d_1$}{
				$m \gets \min\{i:  d_i\le s\}$; \tcp{$1\le m\le n-s+1$} \
				$\mathbf{d_L} \gets (d_{m}\ge d_{m+1} \ge \cdots \ge d_n)$\;
				Remove a copy of $d$ from $\mathbf{d_L}$ if it appears at least once in $\mathbf{d_L}$\;
				\uIf{$\mathbf{d_L}$ has at least $s$ elements}{
					Form all possible $s$ combinations $\mathbf{s_L}$ of $\mathbf{d_L}$, add an element $d'$ with $1\le d'< d$
					into $\mathbf{s_L}$ to form a sequence $\mathbf{s'_L}$ of length $s+1$; Form all possible sequences $\mathbf{s_H}$ of length
					$n-s-1$ consisting of the elements of $\mathbf{d}-(\{d\}\cup \mathbf{s_L})$, add an element $d''=d-d'$ ($1\le d''< d$)
					into $\mathbf{s_H}$ to form a sequence $\mathbf{s'_H}$ of length $n-s$; If both $\mathbf{s'_H}$ and $\mathbf{s'_L}$
					have even sum and are graphical, \Return{\textit{False}}\;
				}
			}
			
		}
	}
	\Return{\textit{True}}\;
	\caption{improved version pseudo-code to test whether a graphical degree sequence is forcibly biconnected}
	\label{algf2_improved}
\end{algorithm}

%
\begin{figure}[!htbp]
	\centering
	\includegraphics[width=8cm]{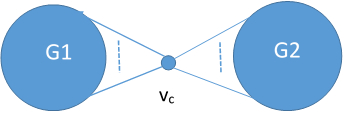}\\
	\caption{A realization $G$ of $\mathbf{d}$ with a cut vertex $v_c$.}
	\label{fig:r1}
\end{figure}
Now we show why Algorithm \ref{algf2_improved} correctly identifies whether the input $\mathbf{d}$ is forcibly biconnected or not.
The conditional test on line 1 is the same as in Algorithm \ref{algf2_basic} since potentially biconnectedness and forcibly
connectedness are necessary conditions for forcibly biconnectedness. The following discussion about the pseudo code
of Algorithm \ref{algf2_improved} from line 3 on all safely assumes that the input $\mathbf{d}$ is potentially biconnected
and forcibly connected.

The conditional test on line 3 works because $d_2+d_n< n$ is a necessary condition for the existence of a cut degree.

\textit{Proof:} Assume there is a cut vertex $v_c$ in a realization $G$ of $\mathbf{d}$. There are two possible situations with
respect to the degree deg($v_c$) of $v_c$.
\begin{itemize}
	\item If deg($v_c$) = $d_1$, then consider two vertices $v_2$ and $v_x$ in the two disjoint induced subgraphs $G_1$ and $G_2$
	of $G$ separated by $v_c$ respectively such that deg($v_2$) = $d_2$ and deg($v_x$) is arbitrary. ($G-\{v_c\}$ is the disjoint union
	of $G_1$ and $G_2$. See Figure \ref{fig:r1}. Here the vertex degrees are all with respect to $G$.) Note that there are
	two possibilities 	regarding which induced subgraph $v_2$ and $v_x$ belong: (1) $v_2\in V(G_1)$,
	$v_x\in V(G_2)$; (2) $v_2\in V(G_2)$, 	$v_x\in V(G_1)$. It is irrelevant to us here which of the two occurs.
	Since $v_2$ and $v_x$ share at most one neighbor in $G$ (which is $v_c$) 	we have deg($v_2$)+deg($v_x$)$-1\le n-2$.
	Thus $d_2+$deg($v_x$)$\le n-1$. Because $d_n\le$ deg($v_x$), we have $d_2+d_n\le n-1$.
	\item If deg($v_c$) $<d_1$, then consider two vertices $v_1$ and $v_x$ in the two induced subgraphs $G_1$ and $G_2$ of $G$
	separated by $v_c$ respectively such that deg($v_1$) = $d_1$ and deg($v_x$) is arbitrary. By the same argument as above we have
	deg($v_1$)+deg($v_x$)$-1\le n-2$. Thus $d_1+$ deg($v_x$)$\le n-1$. Then we have $d_2+d_n\le n-1$ since $d_2\le d_1$ and $d_n\le$ deg($v_x$).
\end{itemize}

Lines 5 to 8 work as follows. If there is a cut vertex $v_c$ in a realization $G$ of $\mathbf{d}$ whose removal results in two
disjoint induced subgraphs $G_1$ and $G_2$, then there are some restrictions on the order of the smaller of these two subgraphs
given $\mathbf{d}$. Without loss of generality we assume from now on $|V(G_2)|\le |V(G_1)|$ so that $G_2$ is
always the smaller of the two disjoint subgraphs.
This piece of code is trying to find out the set $\mathbf{S}$ of potential values of $|V(G_2)|$. Notice that $G_2$ must contain
at least $d_n$ vertices because if it contains at most $d_n-1$ vertices, then each of its vertices will have a degree at most $d_n-1$
in $G$, which contradicts the fact that $d_n$ is the smallest degree. Obviously $G_2$ cannot contain more than $\lfloor (n-1)/2 \rfloor$
vertices since $G_2$ is the smaller of the two induced subgraphs separated by $v_c$. This explains the \textbf{for} loop
lower bound ($d_n$) and upper bound ($\lfloor (n-1)/2 \rfloor$) on line 6. For an integer $s$ between $d_n$ and $\lfloor (n-1)/2 \rfloor$
to be a potential order of $G_2$ (i.e. $|V(G_2)|=s$) we must have $d_{n-s+1}\le s$ since all the vertices in $G_2$ have degrees
at most $s$ in $G$, which makes it necessary for $\mathbf{d}$ to have at least $s$ terms that are at most $s$.
The $s$ smallest terms of $\mathbf{d}$ are $d_{n-s+1}\ge d_{n-s+2} \ge \cdots \ge d_n$.
Thus we must have $d_{n-s+1}\le s$. If $|V(G_2)|=s$ then $|V(G_1)|=n-s-1$. Therefore all the vertices in $G_1$ have degrees
at most $n-s-1$ in $G$. Now other than the requirement that $\mathbf{d}$ has at least $s$ terms that are $\le s$, the other necessary
condition for $\mathbf{d}$ to contain a cut degree is the existence of an \textit{additional} $n-s-1$ terms that are at most $n-s-1$.
Since $s\le n-s-1$ based on the chosen lower and upper bounds for $s$, we must have $d_2\le n-s-1$ as the other necessary
condition.

The functionality of lines 9 to 10 is now clear. If no potential smaller induced subgraph order $|V(G_2)|$ can be found in the range
from $d_n$ to $\lfloor (n-1)/2 \rfloor$, then definitely no cut degree in $\mathbf{d}$ can possibly exist and $\mathbf{d}$ must be forcibly
biconnected.

After testing the necessary conditions the double \textbf{for} loop from lines 11 to 18 tries to find out if a cut degree $d$ in $\mathbf{d}$
exists such that the smaller subgraph $G_2$ after cut has order $s$ from $\mathbf{S}$.

The conditional test on line 13 is to ensure that the largest degree $d^{*}$ in $\mathbf{d}$
that is not the cut degree $d$ does not exceed the order $n-s-1$ of the larger subgraph $G_1$ after cut, should a cut degree $d$
exists. This is because any vertex in $G_1$ can only have the vertices in $G_1$ or
the cut vertex with degree $d$ as its neighbors in $G$. If the vertex with the degree $d^{*}$
occurs in $G_2$, then $d^{*}$ actually should not exceed the order $s$ of $G_2$. However, since $s\le n-s-1$ due to
the upper bound $\lfloor (n-1)/2 \rfloor$ of $s$ on line 6, it is automatically true that $d^{*}$ does not exceed $n-s-1$
when it does not exceed $s$. Now each $d$ in $\mathbf{d}$
is either $d=d_1$ or $d<d_1$. If $d=d_1$, then we need to ensure that $d_2$, which is the largest degree $d^{*}$ in $\mathbf{d}$
other than $d$, does not exceed $n-s-1$, i.e. $d_2\le n-s-1$. This inequality has already been satisfied due to the conditional test
on line 7. If $d<d_1$, then we need to ensure that $d_1$, which is the largest degree $d^{*}$ in $\mathbf{d}$
other than $d$, does not exceed $n-s-1$, i.e. $d_1\le n-s-1$. This completes our justification of the conditional test on line 13.

Lines 14 to 16 construct the sub sequence $\mathbf{d_L}$ of $\mathbf{d}$ consisting exactly of those elements at most
the order $s$ of the smaller subgraph $G_2$. The motivation for this construction is that the $s$ degrees of the vertices
of $G_2$ must all come from $\mathbf{d_L}$. Note that if the cut degree $d\le s$, then we need to remove one copy of $d$
from $\mathbf{d_L}$ since the cut vertex itself is not considered part of $G_2$. This explains why we have line 16.

The conditional test on line 17 is to ensure that there are enough degrees $\le s$ from $\mathbf{d}$ for $G_2$. The sequence
$\mathbf{d_L}$ constitutes the pool of degrees we can select for the vertex degrees of $G_2$. Since $G_2$ has order $s$
and each of its vertices has degree $\le s$, we need $\mathbf{d_L}$ to have at least $s$ elements.

Line 18 performs exhaustive enumerations to find out if there are two graphical degree sequences $\mathbf{s'_H}$
and $\mathbf{s'_L}$, which are the vertex degree sequences of $G-V(G_2)$ and $G-V(G_1)$ respectively, should a cut
degree $d$ exists such that the smaller subgraph $G_2$ is of order $s$ after cut. The vertex degree sequence $\mathbf{s'_L}$
of $G-V(G_1)$
consists of the $s$ degrees of the vertices of $G_2$ in $G$ together with the degree $d'$ equal to the number of adjacencies
of the cut vertex with $G_2$. The vertex degree sequence $\mathbf{s'_H}$ of $G-V(G_2)$ consists of the $n-s-1$ degrees of
the vertices of $G_1$ in $G$ together with the degree $d''=d-d'$ equal to the number of adjacencies of the cut vertex with $G_1$.
Clearly, one level of exhaustive enumeration comes from selecting $s$ degrees $\mathbf{s_L}$ from $\mathbf{d_L}$.
The other level of exhaustive enumeration comes from choosing $d'$ ($1\le d'< d$), which is the number of adjacencies
of the cut vertex with $G_2$. Notationally we have $\mathbf{s'_L}=\mathbf{s_L}\cup \{d'\}$. Once $\mathbf{s_L}$ from
$\mathbf{d_L}$ and $d'$ ($1\le d'<d$) have been chosen, $\mathbf{s'_L}$ and $\mathbf{s'_H}$ are both determined, with
$\mathbf{s'_H}=(\mathbf{d}-(\{d\}\cup \mathbf{s_L}))\cup \{d''\}$. We only need a linear time algorithm such as from \cite{Ivanyi2013}
to test whether $\mathbf{s'_L}$ and $\mathbf{s'_H}$ are both graphical. If during the double \textbf{for} loop from lines 11 to 18
such a pair of graphical degree sequences $\mathbf{s'_L}$ and $\mathbf{s'_H}$ can be found, then we know that the input
$\mathbf{d}$ is not forcibly biconnected (hence returning \textit{False} on line 18) since we have found a cut degree $d$ together
with the degrees for vertices of $G_1$ and $G_2$ after cut and the number of adjacencies $(d'',d')$ of the cut vertex with $G_1$
and $G_2$ respectively. If no such pair of
graphical degree sequences $\mathbf{s'_L}$ and $\mathbf{s'_H}$ can ever be found on line 18, then we know that the
input $\mathbf{d}$ is forcibly biconnected and we should return \textit{True} on line 19.

One thing to note about Algorithm \ref{algf2_improved} is that the sequences mentioned from line 15 to line 18 are all treated as
multisets so that the set difference and set union operations therein should be implemented as multiset operations. We also note that
Algorithm \ref{algf2_improved} performs a test of whether any graphical sequence is forcibly connected or not
at most once (on line 1), while Algorithm \ref{algf2_basic} possibly performs such a test many times (on lines 1 and 5).
This is a supplementary explanation of why Algorithm \ref{algf2_improved} performs much better than Algorithm \ref{algf2_basic} on
many inputs. More detailed run time analysis is presented later.

\subsection{Extensions of the algorithm}
In this section we briefly discuss how to extend Algorithm \ref{algf2_improved} to perform the additional task of listing all
possible cuts of a potentially biconnected and forcibly connected but non forcibly biconnected graphical degree sequence.
We also show that the idea of Algorithm \ref{algf2_improved} to find a cut degree in determining forcibly biconnectedness
can be extended to test forcibly triconnectedness of a graphical degree sequence and beyond.
\subsubsection{Enumeration of all possible cuttings}
It is easy to see that if the input $\mathbf{d}$ is potentially biconnected and forcibly connected but not forcibly
biconnected, we can enumerate all possible cuttings (the cut degree, the numbers of adjacencies of the cut vertex to $G_1$ and
$G_2$ and the degrees of the vertices of $G_1$ and $G_2$ in $G$ using the notation of Section \ref{subsec:pseudo-code}).
We simply need to report such a cutting on line 18 of Algorithm \ref{algf2_improved}
without returning \textit{False} immediately. Such an enumerative algorithm to find all possible cuttings of the input $\mathbf{d}$
can be useful when we want to explore the possible realizations of $\mathbf{d}$ and their properties.

\subsubsection{Testing forcibly $k$-connectedness of $\mathbf{d}$ when $k\ge 3$}
The idea to find a cut degree in Algorithm \ref{algf2_improved} can be extended to find a pair of cut degrees for the
case of $k=3$ and a triple of cut degrees for the case of $k=4$, and so on. Now we need to consider more complicated
situations about whether those cut vertices with potential cut degrees are themselves adjacent or not,
besides the numbers of adjacencies of each of these cut vertices to $G_1$ and $G_2$ and the potential orders of
$G_1$ and $G_2$. However, we believe a careful implementation should still have better performance than the basic
version to test forcibly $k$-connectedness presented in \cite{Wang2018}.

\section{Complexity analysis}
\label{sec:analysis}
In this section we give a rough and conservative analysis of the worst case run time performance of Algorithms
\ref{algf2_basic} and \ref{algf2_improved}. We
also comment on the computational complexity of the decision problem of testing forcibly biconnectedness of graphical
degree sequences.

Let us first consider the basic Algorithm \ref{algf2_basic}.

The test of potentially biconnectedness on line 1 can be performed in linear time using the Wang and Cleitman
characterization \cite{WangKleitman1973}. The test of forcibly connectedness on line 1 has been analyzed in
\cite{Wang2018} with worst case run time probably exponential (around $poly(n)\binom{n}{n/2}$ where $poly(n)$ is some
polynomial of $n$).

The number of iterations of the body of the \textbf{for} loop from lines 3 to 6 depends on the number of possible choices of $d_i$
and $\mathbf{d_S}$. The number of choices for $d_i$ is $O(n)$. The number of choices of $\mathbf{d_S}$ after $d_i$
is chosen depends on the magnitude of $d_i$. Assume the average size of $d_i$ is $n/2$ for the purpose of simplicity
of our analysis. The number of choices of $\mathbf{d_S}$ is then about $\binom{n-1}{n/2}$ since $\mathbf{d_S}$ is to be
chosen from $\mathbf{d}-\{d_i\}$. Once $d_i$ and $\mathbf{d_S}$
have been chosen, the sequence $\mathbf{d'}$ of length $n-1$ can be constructed and tested whether it is a graphical
sequence in $O(n)$ time. After that the testing of forcibly connectedness of $\mathbf{d'}$ can be done in $poly(n-1)\binom{n-1}{(n-1)/2}$
time. Thus we can see that the worst case run time complexity of the \textbf{for} loop from lines 3 to 6 is roughly
$poly(n)\binom{n}{n/2}^2$ under this simplified analysis.

The total run time complexity of Algorithm \ref{algf2_basic} is dominated by the \textbf{for} loop from lines 3 to 6,
which is about $poly(n)\binom{n}{n/2}^2$.

Let us now turn to the improved Algorithm \ref{algf2_improved}.

The worst case run time of line 1 to test potentially biconnectedness and forcibly connectedness is about $poly(n)\binom{n}{n/2}$
as shown above.

It is easy to see that lines 3 to 10 take linear time. The body of the double \textbf{for} loops from line 13 to 18 will
be iterated at most $O(n^2)$ times since we have $O(n)$ distinct $d$ in $\mathbf{d}$ and less than $n/2$ elements
in $\mathbf{S}$. Lines 14 to 16 clearly take $O(n)$ time. Thus, the total run time excluding the test of forcibly
connectedness on line 1 and the exhaustive enumeration on line 18 is $O(n^3)$.

Now we consider line 18. The length of the sequence $\mathbf{d_L}$ could be up to $n$ and the length $s$ of the
sequence $\mathbf{s_L}$ could be up to $\lfloor (n-1)/2 \rfloor$. Then the maximum number of all possible choices of $\mathbf{s_L}$
from $\mathbf{d_L}$ will be up to $\binom{n}{\lfloor (n-1)/2 \rfloor}$, which is exponential. The number of choices for $d'$ is
$O(n)$. Once the choices for $\mathbf{s_L}$ and $d'$ have been made, the two sequences $\mathbf{s'_L}$ and
$\mathbf{s'_H}$ can be constructed and tested whether they are graphical
in $O(n)$ time. This shows that the worst case total run time of line 18 is about $poly(n)\binom{n}{n/2}$.

The most time consuming parts of Algorithm \ref{algf2_improved} are the testing of forcibly connectedness on line 1
and the exhaustive enumeration on line 18. In the worst case both of these parts can take exponential time. The combined
worst case total run time of Algorithm \ref{algf2_improved} is then about $poly(n)\binom{n}{n/2}$, which can be seen to
be better than the worst case run time of Algorithm \ref{algf2_basic} by a factor of $\binom{n}{n/2}$.
Suitable data structures can be established for those multiset sequences mentioned on line 18
to avoid duplicate exhaustive enumerations.
However, a large number of exhaustive enumerations on line 18 before returning can happen for some inputs based on
our computational experiences. 

As for the computational complexity of the problem of deciding forcibly biconnectedness of arbitrary graphical degree
sequences, it is known to be in co-NP as indicated in \cite{Wang2018}. However we do not know if it is co-NP-hard or
if it is inherently harder than the problem of deciding forcibly connectedness of arbitrary graphical degree sequences.
In general one can be interested in whether the problem of deciding forcibly $k$-connectedness is inherently harder than
the problem of deciding forcibly $(k-1)$-connectedness for any fixed $k$.

\section{Computational results}
\label{sec:experiments}
In this section we will first present some results on the experimental evaluation on the performance
of Algorithm \ref{algf2_improved} on randomly generated
graphical degree sequences. We will then provide some enumerative results about the number of forcibly biconnected graphical
degree sequences of given length and the number of forcibly biconnected graphical partitions of a given even integer. Based on
these enumerative results we will make some conjectures about the relative asymptotic
behavior of related functions and the unimodality of certain associated integer sequences.

Before we continue we remark that Algorithm \ref{algf2_basic} can only satisfactorily process inputs with length below 50 most
of the time unless the input is not potentially biconnected or it can be easily determined to be non forcibly connected.
\subsection{Performance evaluations of Algorithm \ref{algf2_improved}}
\label{sec:performance}
Previously in \cite{Wang2018} we evaluated our algorithm to test forcibly connectedness of graphical degree sequences
using randomly generated long inputs with length up to 10000. In anticipation of the greater challenge to test forcibly
biconnectedness (even though the worst case time complexity of the algorithm in \cite{Wang2018} to test forcibly connectedness
and Algorithm \ref{algf2_improved} here differ by at most a polynomial factor as shown in Section
\ref{sec:analysis}), we decide to evaluate Algorithm \ref{algf2_improved} using randomly generated inputs with length up to 1000.

We adopt a similar evaluation methodology as in \cite{Wang2018}. Choose a constant $p_l$ in the range [0.002,0.49] and
a constant $p_h$ in the range [$p_l+0.01$,0.99] and generate 100 random graphical degree sequences of length $n$ with
largest term around $p_hn$ and smallest term around $p_ln$. For each chosen length $n$, the smallest chosen $p_l$ is
slightly adjusted to make sure that the smallest term in each generated random input graphical degree sequence is at least 2
since any graphical degree sequence with the smallest term 1 is not potentially biconnected so the answer is always \textit{False}
for such inputs. The largest $p_l$ is chosen to be 0.49 because any graphical degree sequence of length $n$ with smallest
term at least $0.5n$ is not only forcibly connected (see \cite{Wang2018}) but also forcibly biconnected (see line 3 of Algorithm
\ref{algf2_improved}). The constant $p_h$ is chosen to be less than 1 so that the largest term in any randomly generated input
graphical degree sequence is at most $n-1$.

We implemented our Algorithm \ref{algf2_improved} using C++ and compiled it using g++ with optimization level -O3.
The experimental evaluations are performed on a common Linux workstation.
We run the code on the randomly generated instances and record the average performance
and note the proportion of them that are forcibly biconnected. Table \ref{tab:chosenphpl500} lists the tested $p_l$ and $p_h$ for
the input length $n=500$. The input lengths that are chosen to be tested are $n=20,30,\cdots,100,200,\cdots,1000$.
The chosen $p_l$ and $p_h$ for different input length $n$ are all similar to the case for $n=500$.

\begin{table}[!htb]
	\centering
	\caption{Chosen $p_l$ and $p_h$ in the experimental performance evaluation of Algorithm \ref{algf2_improved} for input length $n=500$.}
	\begin{tabular}{|c|c|}
		\hline
		$p_l$ & $p_h$ \\
		\hline
		\hline
		0.01 & 0.02,0.03,...,0.1,0.2,0.3,...,0.7,\textbf{0.8},\textbf{0.9},\textbf{0.95},\textbf{0.97},\textbf{0.99} \\ \hline
		0.03 & 0.04,0.05,...,0.1,0.2,0.3,...,0.72,\textbf{0.73},\textbf{0.74},\textbf{0.75},\textbf{0.8},\textbf{0.85},\textbf{0.9},\textbf{0.94},0.96,0.99 \\ \hline
		0.06 & 0.08,0.1,0.2,...,0.7,\textbf{0.71},\textbf{0.73},\textbf{0.77},\textbf{0.8},\textbf{0.85},\textbf{0.88},0.89,0.95,0.99 \\ \hline
		0.1 & 0.15,0.25,0.35,...,0.65,\textbf{0.68},\textbf{0.7},\textbf{0.75},\textbf{0.8},0.82,0.85,0.9 \\ \hline
		0.2 & 0.3,0.4,0.5,0.6,\textbf{0.62},\textbf{0.65},\textbf{0.68},0.7,0.8,0.9 \\ \hline
		0.3 & 0.35,0.45,0.55,\textbf{0.56},\textbf{0.58},\textbf{0.6},0.62,0.7,0.8,0.9 \\ \hline
		0.4 & 0.44,0.48,0.52,\textbf{0.53},\textbf{0.54},\textbf{0.55},0.56,0.6,0.7,0.8,0.9 \\ \hline
		0.49 & 0.491,...,0.499,\textbf{0.5},\textbf{0.501},\textbf{0.502},\textbf{0.504},\textbf{0.506},0.508,0.51,0.52,...,0.6,0.7,0.8,0.9 \\ \hline
	\end{tabular}
	\label{tab:chosenphpl500}
\end{table}
We summarize our experimental evaluation results as follows.

1. For all input length $n$ below 100, all instances can be decided instantly (run time $<0.1s$).

2. Starting from $n=100$ up to $n=1000$,
certain inputs start to cause Algorithm \ref{algf2_improved} to run slowly
(run time from a few seconds to a few hours or time out). The longer the length $n$, the more time these difficult instances might
take and the higher percentage of these difficult instances occupy in all 100 input instances for a given triple of ($n,p_l,p_h$).
In detail, difficult instances could occur when (1) $p_l$ is very small (say from 0.005 to 0.05) and $p_h$ is very large
(say from 0.8 to 0.95) or when (2) $p_l$ is slightly below 0.5 and $p_h$ is slightly above 0.5. We note based on our
observations that the former situation (1) often (but not always) corresponds to the cases where forcibly connectedness itself
is hard to decide while the latter situation (2) corresponds to the cases where forcibly connectedness itself is not hard to decide but
forcibly biconnectedness is after confirmation of forcibly connectedness.

3. For each fixed $p_l$, there is a transition range of $p_h$ (shown in \textbf{bold} in Table \ref{tab:chosenphpl500})
such that (1) if $p_h$ is below this transition range, almost all input instances are non forcibly
biconnected; (2) if $p_h$ is above this transition range, almost all input instances are forcibly
biconnected; (3) when $p_h$ increases in the transition range, the proportion of input instances that are forcibly biconnected
grows approximately from 0 to 1. For example, for
input length $n=500$ and $p_l=0.1$, the transition range of $p_h$ is from 0.68 to 0.8. When $p_h$ is below 0.68 almost all input
instances are non forcibly biconnected. When $p_h$ is above 0.8 almost all input instances are forcibly biconnected.
As noted in \cite{Wang2018}, these results indicate relative frequencies, instead of absolute law, of forcibly biconnected graphical
degree sequences among the pool of all graphical degree sequences that can be chosen to test.

To sum up, our implementation of Algorithm \ref{algf2_improved} runs fast on the majority of the tested random inputs with length
up to 1000. Considering that all those inputs with smallest term exactly 1 or at least $n/2$ are excluded from test, which are estimated
to occupy at least 35\% of all possible instances, we are confident that Algorithm \ref{algf2_improved} is efficient on average.
Certain inputs do cause it to perform poorly. In particular, if it is hard to decide whether a potentially biconnected
input is forcibly connected
using the algorithm in \cite{Wang2018}, it is also hard to decide whether it is forcibly biconnected using Algorithm \ref{algf2_improved}
since testing forcibly connectedness is part of testing forcibly biconnectedness. We also remark from our experimental observations that
there are difficult input instances whose hardness come mainly from deciding forcibly connectedness.
That is, the algorithm may take a long time to decide
forcibly connectedness. However, once the algorithm knows the input is potentially biconnected and forcibly connected it can
almost instantly finish the testing of forcibly biconnectedness. There are also difficult input instances that are easy to be determined
to be potentially biconnected and forcibly connected but hard to be further determined to be forcibly biconnected.
Finally, there are some difficult input instances that are hard to be determined to be forcibly connected and, after confirmation
of forcibly connectedness, hard to be further determined to be forcibly biconnected.

\begin{table}[!htb]
	\centering
	\caption{Terminology used in Section \ref{subsec:enum}}
	\begin{tabular}{||c|l||}
		\hline\hline
		Term & Meaning\\
		\hline\hline
		$D(n)$ & number of zero-free graphical sequences of length $n$\\
		\hline
		$D_{k\_c}(n)$ & number of potentially $k$-connected graphical sequences of length $n$ \\
		\hline
		$D_{k\_f}(n)$ & number of forcibly $k$-connected graphical sequences of length $n$ \\
		\hline
		$C_k[n,N]$ & number of potentially $k$-connected graphical degree sequences \\
		& of length $n$ with degree sum $N$ \\
		\hline
		$F_k[n,N]$ & number of forcibly $k$-connected graphical degree sequences \\
		& of length $n$ with degree sum $N$ \\
		\hline
		$L_{k}[n,j]$ & number of forcibly $k$-connected graphical degree sequences of\\
		& length $n$ with largest term $j$ \\
		\hline
		$M_k(n)$ & minimum largest term in any forcibly $k$-connected graphical \\
		& sequence of length $n$ \\
		\hline
		$g(n)$ & number of graphical partitions of even $n$ \\
		\hline
		$g_{k\_c}(n)$ & number of potentially $k$-connected graphical partitions of even $n$ \\
		\hline
		$g_{k\_f}(n)$ & number of forcibly $k$-connected graphical partitions of even $n$ \\
		\hline
		$c_k[n,j]$ & number of potentially $k$-connected graphical partitions of \\
		& even $n$ with $j$ parts \\
		\hline
		$f_k[n,j]$ & number of forcibly $k$-connected graphical partitions of \\
		& even $n$ with $j$ parts \\
		\hline
		$l_k[n,j]$ & number of forcibly $k$-connected graphical partitions of $n$ with largest term $j$\\
		\hline
		$m_k(n)$ & minimum largest term of forcibly $k$-connected graphical partitions of $n$\\
		\hline\hline
	\end{tabular}
	\label{tbl:definitions}
\end{table}

\subsection{Enumerative results}
\label{subsec:enum}
In this section we present some enumerative results related to forcibly biconnected graphical degree sequences
of given length and forcibly biconnected graphical partitions of given even integer. We also make some conjectures
based on these enumerative results. For the reader's convenience, we list the notations and their meanings
used in this section in Table \ref{tbl:definitions}. Note that $C_k[n,N]=c_k[N,n]$ and $F_k[n,N]=f_k[N,n]$ by definition.

In a previous manuscript \cite{Wang2016} we presented efficient algorithms to compute $D(n)$ and $D_{k\_c}(n)$ for fixed $k$ and
have shown that $\lim_{n \to \infty} \frac{D_{1\_c}(n)}{D(n)} = 1$ and
$\lim_{n \to \infty} \frac{D_{2\_c}(n)}{D(n)} \neq 1$. Recently in \cite{Wang2018} we presented enumerative results about
the number $D_{1\_f}(n)$ of forcibly connected graphical degree sequences of length $n$ and conjectured that
$\lim_{n \to \infty} \frac{D_{1\_f}(n)}{D(n)} = 1$. Using the sandwich theorem of calculus we can easily see that
$\lim_{n \to \infty} \frac{D_{2\_f}(n)}{D(n)} \neq 1$ since $D_{2\_f}(n)\leq D_{2\_c}(n)\leq D(n)$.
We adapted the algorithm of Ruskey et al \cite{Ruskey1994} that enumerates
zero-free graphical degree sequences of length $n$ by incorporating the test in Algorithm \ref{algf2_improved} to
enumerate those that are forcibly biconnected. The results together with the proportion of them in all forcibly connected graphical degree
sequences of length $n$ and in all zero-free graphical degree sequences of length $n$ are listed in Table \ref{tab:enumDf2(n)}. From Table
\ref{tab:enumDf2(n)} it looks likely that $D_{2\_f}(n)/D_{1\_f}(n)$ and $D_{2\_f}(n)/D(n)$ both tend to some constant less than 1.
Note that if the conjecture $\lim_{n \to \infty} \frac{D_{1\_f}(n)}{D(n)} = 1$ is true, then
$\lim_{n \to \infty} \frac{D_{2\_f}(n)}{D_{1\_f}(n)} \neq 1$ since we know $\lim_{n \to \infty} \frac{D_{2\_f}(n)}{D(n)} \neq 1$.
\begin{table}[!htb]
	\centering
	\caption{Number of forcibly biconnected graphical sequences of length $n$ and their proportions in forcibly connected
		and zero-free graphical sequences of length $n$ respectively.}
	\begin{tabular}[htbp]{|c||c|c|c|c|}
		\hline
		$n$ & $D_{2\_f}(n)$ & $D_{1\_f}(n)$ & $D_{2\_f}(n)/D_{1\_f}(n)$ & $D_{2\_f}(n)/D(n)$ \\
		\hline
		\hline
		4 & 3  & 6 & 0.500000 & 0.428571 \\ \hline
		5 & 9 & 18 & 0.500000 & 0.450000 \\ \hline
		6 & 30 & 63 & 0.476190 & 0.422535 \\ \hline
		7 & 105 & 216 & 0.486111 & 0.437500 \\ \hline
		8 & 381 & 783 & 0.486590 & 0.437428 \\ \hline
		9 & 1412 & 2843 & 0.496658 & 0.448539 \\ \hline
		10 & 5296 & 10535 & 0.502705 & 0.454397 \\ \hline
		11 & 20010 & 39232 & 0.510043 & 0.461783 \\ \hline
		12 & 76045 & 147457 & 0.515710 & 0.467196 \\ \hline
		13 & 290142 & 556859 & 0.521033 & 0.472392 \\ \hline
		14 & 1110847 & 2113982 & 0.525476 & 0.476649 \\ \hline
		15 & 4264563 & 8054923 & 0.529436 & 0.480473 \\ \hline
		16 & 16411152 & 30799063 & 0.532846 & 0.483750 \\ \hline
		17 & 63284616 & 118098443 & 0.535863 & 0.486662 \\ \hline
		18 & 244489774 & 454006818 & 0.538516 & 0.489220 \\ \hline
		19 & 946101866 & 1749201100 & 0.540877 & 0.491504 \\ \hline
		20 & 3666602417 & 6752721263 & 0.542981 & 0.493542 \\ \hline
		21 & 14229131559 & 26114628694 & 0.544872 & 0.495376 \\ \hline
		22 & 55288167003 & 101153550972 & 0.546577 & 0.497033 \\ \hline
		23 & 215070591363 & 392377497401 & 0.548122 & 0.498537 \\ \hline
		24 & 837503686065 & 1524043284254 & 0.549527 & 0.499908 \\ \hline
		25 & 3264489341370 & 5926683351876 & 0.550812 & 0.501162 \\ \hline
	\end{tabular}
	\label{tab:enumDf2(n)}
\end{table}

In Table \ref{tab:enumD2fdegsums(7)} we show itemized potentially and forcibly biconnected graphical degree sequences
of length 7 based on the degree sum $N$. The counts for $N<14$ are not shown because those counts are all 0. The minimum
$N$ such that $C_2[7,N]$ is nonzero can be easily determined to be 14 using the Wang and Cleitman \cite{WangKleitman1973}
characterization. However, the minimum $N$ such that $F_2[7,N]$ is nonzero (which is equal to 20, not 14) is not obvious.
In general for given $n$ if the minimum $N$ such
that $F_2[n,N]$ is nonzero can be easily determined, it can be added into Algorithm \ref{algf2_improved} so that for any input
$\mathbf{d}$ with the sum of its terms less than this minimum the algorithm can immediately return \textit{False}.
The largest degree sum is 42 for any graphical degree sequence of length 7. Notice that both the nonzero $C_2[7,N]$ values
and the nonzero $F_2[7,N]$ values when $N$ increases form a unimodal sequence. Also notice that $C_2[7,N]=F_2[7,N]$
when $28\le N\le 42$. This shows that for even $N$ between 28 and 42, all potentially biconnected graphical partitions of $N$
with 7 parts are also forcibly biconnected. That is, 28 is the smallest $N$ such that all potentially biconnected graphical partitions
of $N$ with 7 parts are also forcibly biconnected.
\begin{table}[!htb]
	\centering
	\caption{Number of potentially (row $C_2[7,N]$) and forcibly (row $F_2[7,N]$) biconnected graphical degree sequences of length 7 with given degree sum $N$.}
	\begin{tabular}[htbp]{|c||c|c|c|c|c|c|c|c|c|c|c|c|c|c|c|}
		\hline
		degree sum $N$ & 14 & 16 & 18 & 20 & 22 & 24 & 26 & 28 & 30 & 32 & 34 & 36 & 38 & 40 & 42 \\ \hline
		\hline
		$C_2[7,N]$ & 1 & 1 & 3 & 7 & 14 & 17 & 18 & 19 & 16 & 12 & 8 & 5 & 2 & 1 & 1 \\
		\hline
		$F_2[7,N]$ & 0 & 0 & 0 & 2 & 8 & 14 & 17 & 19 & 16 & 12 & 8 & 5 & 2 & 1 & 1 \\ \hline
	\end{tabular}
	\label{tab:enumD2fdegsums(7)}
\end{table}

In Table \ref{tab:enumD2flargestpart(15)} we show itemized numbers of forcibly biconnected graphical degree
sequences of length 15 based on the largest degree. The counts for largest degrees less than 7 are not shown because those
counts are all 0. Any graphical degree sequence of length 15 has a largest term at most 14. From the table we can see that the counts
decrease when the largest degree decreases. For other degree sequence lengths from 5 to 25 we observed
similar behavior. The table also indicates that there are no forcibly biconnected graphical
degree sequences of length 15 with largest degree less than 7. In fact, we can define $M_2(n)$ to be the minimum largest term
in any forcibly biconnected graphical sequence of length $n$. This is, $M_2(n)\doteq$ min\{$\Delta$: $\Delta$ is the largest term
of some forcibly biconnected graphical degree sequence of length $n$\}. Clearly we have $M_2(n)\le n/2$ since for even $n$
the sequence $n/2,n/2,\cdots,n/2$ of length $n$ is forcibly biconnected. In \cite{Wang2018}  we have defined $M(n)$ to be the
minimum largest term in any forcibly connected graphical sequence of length $n$ ($M(n)=M_1(n)$ based on the notation
$M_k(n)$ in Table \ref{tbl:definitions}). By definition we obviously have
$M_2(n)\ge M(n)$ so that we also have $M_2(n)>c\sqrt{n}$ for all sufficiently large $n$ and some constant $c>0$ based on a
lower bound of $M(n)$ shown in \cite{Wang2018}. If $M_2(n)$ can be easily calculated, it can be added into Algorithm
\ref{algf2_improved} so that any input $\mathbf{d}$ of length $n$ with largest term less than $M_2(n)$ can be immediately
decided to be non forcibly biconnected.
We show the values of $M_2(n)$ and $M(n)$ based on our enumerative results in Table \ref{tab:M2n}.
\begin{table}[!htb]
	\centering
	\caption{Number $L_2[15,j]$ of forcibly biconnected graphical degree sequences of length 15 with given largest term $j$.}
	\begin{tabular}[htbp]{|c||c|c|c|c|c|c|c|c|}
		\hline
		largest part $j$ & 14 & 13 & 12 & 11 & 10 & 9 & 8 & 7  \\ \hline
		\hline
		$L_2[15,j]$ & 2113982 & 1335151 & 573980 & 185510 & 45951 & 8689 & 1202 & 98 \\ \hline
	\end{tabular}
	\label{tab:enumD2flargestpart(15)}
\end{table}

\begin{table}[!htb]
	\centering
	\caption{Minimum largest term $M_2(n)$(resp. $M(n)$) of forcibly biconnected (resp. connected) graphical sequences of length $n$.}
	\begin{tabular}[htbp]{|c||c|c|c|c|c|c|c|c|c|c|c|}
		\hline
		$n$ & 4 & 5 & 6 & 7 & 8 & 9 & 10 & 11 & 12 & 13 & 14 \\
		\hline
		$M_2(n)$ & 2 & 2 & 3 & 3 & 4 & 4 & 4 & 5 & 5 & 6 & 6 \\
		\hline
		$M(n)$ & 2 & 2 & 3 & 3 & 3 & 4 & 4 & 5 & 5 & 5 & 6 \\ \hline \hline
		$n$ & 15 & 16 & 17 & 18 & 19 & 20 & 21 & 22 & 23 & 24 & 25 \\
		\hline
		$M_2(n)$ & 7 & 7 & 7 & 8 & 8 & 9 & 9 & 10 & 10 & 10 & 11 \\
		\hline
		$M(n)$ & 6 & 6 & 7 & 7 & 7 & 7 & 8 & 8 & 8 & 8 & 8 \\ \hline
	\end{tabular}
	\label{tab:M2n}
\end{table}

We also incorporated our Algorithm \ref{algf2_improved} into the highly efficient Constant Amortized Time (CAT) algorithm of
Barnes and Savage \cite{BarnesSavage1997} to generate forcibly biconnected graphical partitions of even $n$. The results for $n$ up to
170 together with the proportion of them in all forcibly connected graphical partitions of $n$ are listed in Table \ref{tab:enumgf2(n)}.
For the purpose of saving space we only show the results in increments of 10 for $n$. From the table it seems
reasonable to conclude that the proportion $g_{2\_f}(n)/g_{1\_f}(n)$ will decrease when $n$ is beyond some small threshold
and it might tend to the limit 0. Previously in \cite{Wang2018} we have conjectured that $\lim_{n \to \infty} \frac{g_{1\_f}(n)}{g(n)}=0$.
With the trend of $g_{2\_f}(n)/g_{1\_f}(n)$ in Table \ref{tab:enumgf2(n)}, we are almost certain that $\lim_{n \to \infty} \frac{g_{2\_f}(n)}{g(n)}=0$.

\begin{table}[!htb]
	\centering
	\caption{Number of forcibly biconnected graphical partitions of $n$ and their proportions in all forcibly connected graphical partitions of $n$.}
	\begin{tabular}[htbp]{|c||c|c|c|}
		\hline
		$n$ & $g_{2\_f}(n)$ & $g_{1\_f}(n)$ & $g_{2\_f}(n)/g_{1\_f}(n)$ \\
		\hline
		\hline
		10 & 2 & 8 & 0.250000 \\ \hline
		20 & 10 & 81 & 0.123457 \\ \hline
		30 & 55 & 586 & 0.093857 \\ \hline
		40 & 262 & 3308 & 0.079202 \\ \hline
		50 & 1062 & 15748 & 0.067437 \\ \hline
		60 & 4171 & 66843 & 0.062400 \\ \hline
		70 & 14445 & 256347 & 0.056349 \\ \hline
		80 & 47586 & 909945 & 0.052295 \\ \hline
		90 & 147132 & 3026907 & 0.048608 \\ \hline
		100 & 430709 & 9512939 & 0.045276 \\ \hline
		110 & 1217258 & 28504221 & 0.042704 \\ \hline
		120 & 3285793 & 81823499 & 0.040157 \\ \hline
		130 & 8621222 & 226224550 & 0.038109 \\ \hline
		140 & 21874986 & 604601758 & 0.036181 \\ \hline
		150 & 54077294 & 1567370784 & 0.034502 \\ \hline
		160 & 130279782 & 3951974440 & 0.032966 \\ \hline
		170 & 306808321 & 9714690421 & 0.031582 \\ \hline
	\end{tabular}
	\label{tab:enumgf2(n)}
\end{table}

When generating all forcibly biconnected graphical partitions of $n$ we can also output the itemized counts based on the
number of parts or the largest part. In Table \ref{tab:enumgf2parts(30)} we show itemized counts of potentially biconnected
and forcibly biconnected graphical partitions of 30 based on the number of parts. The counts $c_2[30,j]$ and $f_2[30,j]$
for which the number of parts $j$ is less than
6 or greater than 15 are not shown because those counts are all 0. When $n$ is large, the minimum number of parts $j$
for which $c_2[n,j]$ and $f_2[n,j]$ are both nonzero is clearly the smallest positive integer $j(n)$ such that
$j(n)(j(n)-1)\ge n$. The largest number of parts $j$ for which $c_2[n,j]$ is nonzero is clearly $n/2$ since the sequence $2,\cdots,2$
($n/2$ copies) is potentially biconnected. The largest number of parts $j$ for which $f_2[n,j]$ is nonzero does not
appear to be easily computable. Note that the nonzero values of $c_2[30,j]$ and $f_2[30,j]$ both form a unimodal sequence
when $j$ increases. Also note that $c_2[30,j]=f_2[30,j]$ for $j=6,7$. This shows that all potentially biconnected graphical
partitions of 30 with 6 or 7 parts are also forcibly biconnected. That is, 7 is the largest number of parts $j$ such that all potentially
biconnected graphical partitions of 30 with $j$ parts are also forcibly biconnected.
The fact that $c_2[30,7]=f_2[30,7]=16$ agrees with the result
that $C_2[7,30]=F_2[7,30]=16$ in Table \ref{tab:enumD2fdegsums(7)}.

\begin{table}[!htb]
	\centering
	\caption{Number of potentially (row $c_2[30,j]$) and forcibly (row $f_2[30,j]$) biconnected graphical partitions of 30 with given number of parts $j$.}
	\begin{tabular}[htbp]{|c||c|c|c|c|c|c|c|c|c|c|}
		\hline
		number of parts $j$  & 6 & 7 & 8 & 9 & 10 & 11 & 12 & 13 & 14 & 15 \\ \hline
		\hline
		$c_2[30,j]$ & 1 & 16 & 44 & 54 & 30 & 15 & 7 & 3 & 1 & 1 \\
		\hline
		$f_2[30,j]$ & 1 & 16 & 30 & 8 & 0 & 0 & 0 & 0 & 0 & 0 \\ \hline
	\end{tabular}
	\label{tab:enumgf2parts(30)}
\end{table}

\begin{table}[!htb]
	\centering
	\caption{Number $l_2[30,j]$ of forcibly biconnected graphical partitions of 30 with given largest term $j$.}
	\begin{tabular}[htbp]{|c||c|c|c|c|c|}
		\hline
		largest part $j$ & 4 & 5 & 6 & 7 & 8 \\ \hline
		\hline
		$l_2[30,j]$ & 2 & 13 & 23 & 13 & 4 \\ \hline
	\end{tabular}
	\label{tab:enumgf2largestpart(30)}
\end{table}

In Table \ref{tab:enumgf2largestpart(30)} we show itemized counts of forcibly biconnected graphical partitions of 30 based
on the largest term. The counts $l_2[30,j]$ for which the largest part $j$ is less than 4 or greater than 8 are not shown
since those counts are all 0. The nonzero counts $l_2[30,j]$ form a unimodal sequence. Given $n$ the largest $j$ for which $l_2[n,j]$
is nonzero can be easily determined by solving the inequality $n\ge 2(j+1)-4+2j$ using the Wang and Cleitman
\cite{WangKleitman1973} characterization since a potentially biconnected graphical degree sequence of length $j+1$
and largest term $j$ is forcibly biconnected. However,
the smallest $j$ for which $l_2[n,j]$ is nonzero does not appear to be easily computable.
If we use $m_2(n)$ to denote this smallest $j$, i.e. the minimum largest term of any forcibly biconnected
graphical partitions of $n$, then we clearly have $2\le m_2(n) \le \sqrt{n}$. We show the values of $m_2(n)$ for some small
$n$ based on our enumerative results in Table \ref{tab:m2n}. The fact that $m_2(30)=4$ agrees with the fact that in
Table \ref{tab:enumgf2largestpart(30)} the counts $l_2[30,j]$ are all 0 for largest part $j$ less than 4. In \cite{Wang2018}
we have defined $m(n)$ as the minimum largest term of any forcibly connected graphical partition of $n$ ($m(n)=m_1(n)$
in the notation $m_k(n)$ of Table \ref{tbl:definitions}). Clearly we have $m(n)\le m_2(n)$ by definition. We can see that $m_2(n)$
happen to agree with $m(n)$ for those $n$ listed in Table \ref{tab:m2n}.
\begin{table}[!htb]
	\centering
	\caption{Minimum largest term $m_2(n)$ (resp. $m(n)$) of forcibly biconnected (resp. connected) graphical partitions of $n$.}
	\begin{tabular}[htbp]{|c||c|c|c|c|c|c|c|c|c|c|}
		\hline
		$n$ & 10 & 20 & 30 & 40 & 50 & 60 & 70 & 80 & 90 & 100 \\ \hline
		\hline
		$m_2(n)$ & 2 & 3 & 4 & 5 & 5 & 6 & 6 & 6 & 7 & 7  \\
		\hline
		$m(n)$ & 2 & 3 & 4 & 5 & 5 & 6 & 6 & 6 & 7 & 7  \\ \hline
	\end{tabular}
	\label{tab:m2n}
\end{table}

\subsection{Questions and conjectures}
\label{sec:conjectures}
Based on the obtained enumerative results we ask the following questions and make certain conjectures:

1. What is the growth order of $D_{2\_f}(n)$ relative to $D_{1\_f}(n)$ and $D(n)$? We conjecture that both $\frac{D_{2\_f}(n)}{D_{1\_f}(n)}$
and $\frac{D_{2\_f}(n)}{D(n)}$ tend to a constant less than 1. Is it true that $\lim_{n \to \infty} \frac{D_{2\_f}(n)}{D_{2\_f}(n-1)}
= 4$? Furthermore, we conjecture that $\frac{D_{2\_f}(n)}{D_{1\_f}(n)}$ and $\frac{D_{2\_f}(n)}{D(n)}$ are both monotonically
increasing when $n\ge 8$.
What can be said about the relative orders of $D_{k\_c}(n)$, $D_{k\_f}(n)$ and $D(n)$ when $k\geq 3$?

2. What is the growth order of $g_{2\_f}(2n)$ relative to $g_{1\_f}(2n)$ and $g(2n)$? We conjecture that
$\lim_{n \to \infty} \frac{g_{2\_f}(2n)}{g_{1\_f}(2n)}= 0$ and $\lim_{n \to \infty} \frac{g_{2\_f}(2n)}{g(2n)}= 0$. Furthermore,
we conjecture that $g_{2\_f}(2n)/g_{1\_f}(2n)$ is monotonically decreasing when $n\ge 5$. Is it true that
$\lim_{n \to \infty} \frac{g_{k\_f}(2n)}{g_{k-1\_f}(2n)}= 0$ for fixed $k\ge 3$?

3. We conjecture that the numbers $F_2(n,N)$ of forcibly biconnected graphical degree sequences of length $n$ with degree
sum $N$, when $N$ runs through $2n-2,2n,\cdots,n(n-1)$, give a unimodal sequence. (There may be
some zeros in the sequence.)

4. Let $t(n)$ be the smallest positive integer such that $t(n)(t(n)-1)\ge n$.
We conjecture that the numbers $f_2(n,j)$ of forcibly biconnected graphical partitions
of $n$ with $j$ parts, when $j$ runs through $t(n), t(n)+1, \cdots, n/2$, give a unimodal sequence. (There may be
some zeros in the sequence.)

5. What is the growth order of $M_2(n)$, the minimum largest term in any forcibly biconnected graphical sequence of
length $n$? Is there a constant $C>0$ such that $\lim_{n \to \infty} \frac{M_2(n)}{n} = C$? Is it true that $M_2(n)=\Theta(M(n))$?
Is there an efficient algorithm to compute $M_2(n)$? What can be said about the growth order of $M_k(n)$ for fixed
$k\ge 3$?

6. What is the growth order of $m_2(n)$, the minimum largest term in any forcibly biconnected graphical partition of
$n$? Is there a constant $C>0$ such that $\lim_{n \to \infty} \frac{m_2(n)}{\sqrt{n}} = C$?
Is there an efficient algorithm to compute $m_2(n)$? What can be said about the growth order of $m_k(n)$ for fixed
$k\ge 3$?

7. We conjecture that the numbers $l_2(n,\Delta)$ of forcibly biconnected graphical partitions of an even $n$
with largest part exactly $\Delta$, when $\Delta$ runs through $m_2(n),m_2(n)+1,\cdots,n/2$, give a unimodal sequence.
(There may be some zeros in the sequence.)

\section{Conclusions}
\label{sec:conclusion}
In this paper we presented an efficient algorithm to test whether a given graphical degree sequence is forcibly
biconnected or not. We also discussed the possible extension of the idea used in the algorithm to test forcibly $k$-connectedness of
graphical degree sequences for fixed $k\ge 3$. The worst case run time complexity of the algorithm is
exponential. However, extensive performance evaluations on a wide range of random graphical degree sequences
demonstrate its average case efficiency.
We incorporated this testing algorithm into existing algorithms that enumerate zero-free graphical
degree sequences of length $n$ and graphical partitions of an even integer $n$ to obtain some enumerative results
about the number of forcibly biconnected graphical degree sequences of length $n$ and forcibly biconnected graphical
partitions of $n$. Questions and conjectures based on these enumerative results are proposed which warrant a lot of
further research.

\section{Acknowledgements}
This research has been partially supported by a research seed grant of Georgia Southern University.
\bibliographystyle{plain}
\bibliography{testF2}

\end{document}